\documentclass[11pt]{amsart}
\usepackage{amssymb,latexsym}
\newtheorem{theorem}{Theorem}[section]
\newtheorem{lemma}[theorem]{Lemma}
\newtheorem{proposition}[theorem]{Proposition}

\theoremstyle{definition}

\theoremstyle{remark}
\newtheorem{remark}[theorem]{Remark}

\numberwithin{equation}{section}
\begin{document}
\title[semilinear problem on H-type group]{The Liouville theorem  on H-type groups}

\author{Chuanyang Li$^1$, Juan Zhang$^1$ and Peibiao Zhao$^{1,*}$}


\begin{abstract}
In this paper we obtain a Liouville type theorem to the semilinear subcritical elliptic equation on H-type groups. The semilinear subcritical elliptic equation studied in this paper is a generalization of a classical semilinear subcritical elliptic equation on the Heisenberg group. The proofs are based on an {\it a priori} integral estimate and a generalized differential identity which  found by Jerison and Lee [J. Diff. Geom, 29 (1989)].\\
\end{abstract}
 \maketitle



\section{Introduction}

\setcounter{equation}{0}

In this paper, we consider solutions to the following critical semilinear elliptic equation

\begin{equation}\label{1.1}
   -\Delta_\mathbb{G}u=2n^2u^q \quad \text{in} \quad \Omega,
\end{equation}
where $\Omega$ is a domain in H-type group $\mathbb{G}$ (also known as group of Heisenberg type), and  $u$ is a real, smooth and nonnegative function defined in $\Omega$,
while $\triangle_{\mathbb{G}}u=u_{\alpha\overline{\alpha}}+u_{\overline{\alpha}\alpha}$ is the H-type laplacian of $u$, and $q^*=\frac{Q+2}{Q-2}$with$Q=2n+2m$ the homogeneous dimension of  $\mathbb{G}$.  Our main results are
an entire Liouville type theorem and an integral estimate for solutions to (\ref{1.1}). Indeed we have

\begin{theorem}\label{Thm1}
Let $\Omega=\mathbb{G}$ be the whole space. Then the equation (\ref{1.1}) has no positive  solution,
namely, the only nonnegative solution of (\ref{1.1}) is the trivial one.
\end{theorem}

\begin{theorem}\label{Thm2}
 Let $1<q<q^*$, $B_{4r}(\xi_0)\subset \Omega$ be any ball centered at $\xi_0$ with radio $4r$.
 Then any positive solution $u$ of (\ref{1.1}) satisfies:

\begin{equation}\label{1.2}
  \int_{B_r(\xi_0)}  u^{3q-q^*} \leq C\,r^{Q-2\times[\frac{3q-q^*}{q-1}+\frac{2}{(q-q^*)(n+m-1)+2}]} ,
\end{equation}
with some positive constant $C$ depending only on $n$ ,$m$ and $q$.
\end{theorem}

We can prove theorem \ref{Thm1} through theorem \ref{Thm2}. For  $1<q<q^*$, we see $Q-2\times[\frac{3q-q^*}{q-1}+\frac{2}{(q-q^*)(n+m-1)+2}]<0$. So if $u$ be a positive solution of (\ref{1.1}) with $\Omega=\mathbb{G}$,
taking  $r\rightarrow +\infty$ in  (\ref{1.2}) we have

\begin{equation}\label{1.3}
   \int_{\mathbb{G}}  u^{3q-q^*} \leq 0.
\end{equation}
This contradiction indicates that theorem \ref{Thm1} is valid.

The analogue of the equation (\ref{1.1}) in the Heisenberg Group
\begin{equation}\label{1.4}
  -\Delta_{\mathbb{H}^n}=2n^2u^q
\end{equation}
has been deeply studied in
decades. The (\ref{1.4}) on the Heisenberg Group comes from the CR Yamabe problem on Heisenberg Group. The CR Yamabe problem was initiated and studied preliminaryly by Jerison and Lee on \cite{JL1987}-\cite{JL1989}. In the subcritical case, for $1<q<q_*$ where $q_*=\frac{Q+2}{Q-2}$ and $Q=2n+2$ is the homogeneous dimension, the nonexistence of positive entire solutions (Liouville theorem) and a point wise estimate near an isolated singularity for (\ref{1.4}) was obtained by Ma and Ou \cite{OM2023}. The core of their proof relies on a novel a priori integral estimate, which is derived from a generalized and more transparent version of a key identity originally introduced by Jerison and Lee. This paper will also conduct the proof following their method. For
 \begin{equation}\label{1.5}
  -\Delta_{\mathbb{H}^n}=2n^2u^{q_*},
\end{equation}
there are positive solutions of (\ref{1.5}) given by Jerison and Lee \cite{JL1988}
\begin{equation}\label{1.6}\mathcal{U}_{\lambda,\mu}(z,t):=\frac{C}{|t+i|z|^2+z\cdot\mu+\lambda|^n},\end{equation}
for some $\lambda\in\mathbb{C}$,$\mu\in\mathbb{C}^n$,$\mathrm{Im}(\lambda)>\frac{|\mu|^2}{4}$ and $C=C(n,\lambda)>0$.
Catino and Li \cite{CL2023} also obtained a complete classification of positive solutions for (\ref{1.5}).

The analogue of the equation is (\ref{1.1}) in the Euclidean space
\begin{equation}\label{1.7}-\Delta u=u^{n^*-1}\quad\mathrm{in~}\mathbb{R}^n,\end{equation}
where $n^*=\frac{2n}{n-2}$.
Equation (\ref{1.3}) is connected to the Yamabe problem in Riemannian geometry \cite{JL1987} and the extremals in the Sobolev inequality \cite{Ro2022}.From \cite{Ro1966},\cite{Ta1976} and \cite{Au1976} we know that there are the solutions of (\ref{1.7})
\begin{equation}\label{1.8}\mathcal{V}_{\lambda,x_0}(x):={\left(\frac{\lambda\sqrt{n(n-2)}}{\lambda^2+|x-x_0|^2}\right)^{\frac{n-2}{2}}},\end{equation}
for $\lambda>0$ and $x_0\in\mathbb{R}^n$. For (\ref{1.7}),in \cite{CG1989}Caffarelli, Gidas and Spruck obtained a classification of solutions, namely (\ref{1.8}) are the only positive solutions to (\ref{1.7}). In the subcritical case
\begin{equation}\label{1.9}-\Delta u=u^q\quad\mathrm{~in~}\mathbb{R}^n,\end{equation}
where $1<q<n^*-1$. Gidas and Spruck \cite{GS1981} obtained that the only nonnegative solution of (\ref{1.9}) is the trivial one.

In \cite{GV2001} Garofalo and Vassilev study that symmetry properties and uniqueness of positive entire solutions to the Yamabe-type equation
\begin{equation}\label{1.10}
\mathcal{L}u = -u^{\frac{Q+2}{Q-2}}
\end{equation}
on groups of Heisenberg type, including Iwasawa-type groups. The authors establish that if a non-trivial positive entire solution possesses partial symmetry, then it must actually exhibit cylindrical symmetry, and furthermore, it must belong to the explicit one-parameter family of solutions \( K_\epsilon \) (modulo group translations and dilations). Their work generalizes the celebrated results of Jerison and Lee on the Heisenberg group to the broader setting of groups of Heisenberg type, providing a complete symmetry and uniqueness theory for positive entire solutions of the Yamabe equation in this context.

\section{Preliminaries,Notations And Generalization of Jerison-Lee's identity in H-type group}

\setcounter{equation}{0}
\setcounter{theorem}{0}

We first give a brief introduction to the H-type group  $\mathbb{G}$ and some notation.
Let $U^{(k)}$ is a family of skew-Hermitian matrix. $U^{(k)}$ can be written as $A^{(k)}+iB^{(k)}$. Let
$$V^{(k)}=\begin{pmatrix}
A^{(k)}& -B^{(k)} \\
B^{(k)} & A^{(k)}
\end{pmatrix}$$
And $V^{(k)}$ satisfy the following conditions. $V^{(k)}$ is an antisymmetric orthogonal matrix, $\forall j=1,2,\cdot\cdot\cdot,m$;$V^iV^j+V^jV^i=0, \forall i,j\in\{1,2,\cdots,m\},i\neq j$.
   We consider
   $$\mathbb{G}:=\mathbb{C}^n\times\mathbb{R}^m$$
   with coordinates ($w,\,t$) and group law $\circ$: $$(w,t)\cdot(w^{\prime},t^{\prime})=\left(w+w^{\prime},t+t^{\prime}+2\operatorname{Re}\left(\sum_{k=1}^m\langle w,\overline{U^{(k)}}w^{\prime}\rangle e_k\right)\right),$$
     $$\text{for}\,\,(w,t),\,(w^{\prime},t^{\prime})\in \mathbb{C}^n\times \mathbb{R}^m, $$
\noindent where $\left\langle\cdot,\cdot\right\rangle$ is the inner product of Hermite and $e^k$ is the unit vector. We define, for $\xi=(w,t)\in\mathbb{G}$, the distance from the origin as $\mathrm{d}(\xi)=(|w|^{4}+16|t|^{2})^{\frac{1}{4}}$. We use the notation $B_R(\xi)$ for the metric ball centered at $\xi\in\mathbb{H}^{n}$ with radius $R>0$, that is
$$B_R(\xi)=\{\zeta\in\mathbb{H}^n:d(\xi,\zeta)<R\}$$
It is well-known and it is important to recall that the volume of a metric ball is given by
$$|B_r(\xi)|=Cr^Q,$$
where $C>0$ is a positive constant,$Q=2n+2m$ and $\left|\cdot\right|$ denotes the Lebesgue measure.

We define the following left-invariant vector fields in $\mathbb{G}$
$$Z_\alpha=\frac{\partial}{\partial w_\alpha}+\sum_{k=1}^n\left(\sum_{\beta=1}^nU_{\alpha\beta}^{(k)}{\overline{w}_\beta}\right)\frac{\partial}{\partial t_k},\quad Z_{\overline{\alpha}}=\frac{\partial}{\partial{\overline{w}_\alpha}}+\sum_{k=1}^n\left(\sum_{\beta=1}^n\overline{U_{\alpha\beta}^{(k)}}w_\beta\right)\frac{\partial}{\partial t_k}$$
for $\alpha=1,\ldots,n$ and $k=1,\ldots,m$.

We would indicate the derivatives of functions or vector fields with indices preceded by a comma, to avoid confusion. For a smooth function $f:\mathbb{G}\to\mathbb{C}$ we denote its derivatives by
$$f_{\alpha}=Z_{\alpha}f,\quad f_{\overline{\alpha}}=Z_{\overline{\alpha}}f,\quad f_{k0}=\frac{\partial f}{\partial t_{k}},$$
$$\quad f_{\alpha\overline{\beta}}=Z_{\overline{\beta}}\left(Z_{\alpha}f\right),\quad f_{k0,\alpha}=Z_{\alpha}\left(\frac{\partial f}{\partial t_{k}}\right),$$
and so on.There hold the following commutation rules
$$f_{\alpha\beta}-f_{\beta\alpha}=0,\quad f_{\alpha\overline{\beta}}-f_{\overline{\beta}\alpha}=2U_{\alpha\beta}^{(k)}f_{k0},$$
$$\quad f_{k0,\alpha}-f_{\alpha,k0}=0,
\quad f_{\alpha\beta\overline{\gamma}}-f_{\alpha\overline{\gamma}\beta}=2U_{\beta\gamma}^{(k)}f_{k0,\alpha}.$$
Where here and in the sequel we use the Einstein notation sum for the Greek indices $1\leq\alpha,\beta,\gamma\leq n$ and$1\leq k,k'\leq m$.
Moreover, we define
$$|\partial f|^2:=\sum_{\alpha=1}^nf_\alpha f_{\overline{\alpha}}=f_\alpha f_{\overline{\alpha}}.$$
Given $u>0$ a solution of (\ref{1.1}), we consider the auxiliary function $f$ defined as follows $e^f= u^{\frac{1}{n}}$.
Take $q=q^*+ \frac{p}{n+m-1}$, then $f$ satisfies the following equation
\begin{equation}\label{2.1}
\mathbf{Re} f_{\alpha\overline{\alpha}}=-n|\partial f|^2-ne^\frac{n(2+p)f}{n+m-1}.
\end{equation}
Define the tensors
\begin{equation}\label{2.2}
\begin{split}
 D_{\alpha\beta}            =& f_{\alpha\beta}-2f_{\alpha}f_{\beta}, \qquad\qquad\qquad  D_{\alpha}=D_{\alpha\beta}f_{\overline{\beta}},\\
 E_{\alpha\overline{\beta}} =& f_{\alpha\overline{\beta}}-\frac{1}{n}f_{\gamma\overline{\gamma}}\delta_{\alpha\overline{\beta}},
                              \qquad\quad \qquad E_{\alpha}=E_{\alpha\overline{\beta}}f_{\beta}, \\
 G_\alpha=&\frac{U_{\beta\beta}^{(k)}f_{k0,\alpha}}{n}-\frac{U_{\beta\beta}^{(k)}f_{k0}}{n}f_\alpha+e^{(2+p)f}f_\alpha+\mid\partial f\mid^2f_\alpha.
\end{split}
\end{equation}
Denote the function
$$g=\mid\partial f\mid^2+e^{\frac{n}{n+m-1}(2+p)f}-\frac{U_{\alpha\alpha}^{(k)}f_{k0}}{n}.$$Then we can rewrite the equation (\ref{2.1}) as
\begin{equation}\label{2.3}
 f_{\alpha\overline{\alpha}}=-ng.
\end{equation}
\noindent Moreover, we  notice that
\begin{equation}\label{2.4}
\begin{split}
\,&  E_{\alpha\overline{\beta}} = f_{\alpha\overline{\beta}}+g\delta_{\alpha\overline{\beta}},\qquad\,\,\qquad
     E_{\alpha}= f_{\alpha\overline{\beta}}f_{\beta}+gf_{\alpha},\\
\,&  D_{\alpha}=f_{\alpha\beta}f_{\overline{\beta}}-2|\partial f|^2f_{\alpha},\qquad
     G_\alpha=\frac{U_{\beta\beta}^{(k)}f_{k0,\alpha}}{n}+gf_\alpha.
\end{split}
\end{equation}
And by observation we find$f_{\alpha}E_{\overline{\alpha}}=f_{\overline{\alpha}}E_{\alpha} $. This means it's real.
And by Straight calculations, we have
\begin{equation}\label{2.5}
 (|\partial f|^2)_{,\overline{\alpha}} =D_{\overline{\alpha}}+E_{\overline{\alpha}}+\overline{g}f_{\overline{\alpha}}-2f_{\overline{\alpha}}e^\frac{n(2+p)f}{n+m-1},
\end{equation}
and
\begin{equation}\label{2.6}
\begin{split}
g_{k0}&=-\frac{n}{U_{\gamma\gamma}^{(k)}}f_{\alpha}G_{\overline{\alpha}}+\frac{n}{U_{\gamma\gamma}^{(k)}}f_{\overline{\alpha}}G_{\alpha}+2f_{k0}\mid\partial f\mid^2\\&\quad+\frac{n}{n+m-1}(2+p)f_{k0}e^{\frac{n}{n+m-1}(2+p)f}-\frac{U_{\gamma\gamma}^{(k)}f_{k0,k'0}}{n}.
\end{split}
\end{equation}
And by \label{2.5} we find
\begin{equation}\label{2.7}
\begin{split}
g_{\overline{\alpha}} =& D_{\overline{\alpha}}+E_{\overline{\alpha}}+G_{\overline{\alpha}}+\frac{np-2m+2}{n+m-1}f_{\overline{\alpha}}e^\frac{n(2+p)f}{n+m-1}.
\end{split}
\end{equation}

By generalizing the equations on the Heisenberg group based on the above observations, the following equality is obtained
\begin{proposition}\label{Pro-1}
\begin{equation}\label{2.8}
\begin{split}
\begin{aligned}
&\mathcal{M}=\mathbf{Re}\{Z_{\overline{\alpha}}\{e^{2(n-1)f}[(g+3\frac{U_{\gamma\gamma}^{(k)}f_{k0}}{n})E_{\alpha}+\left(g-\frac{U_{\gamma\gamma}^{(k)}f_{k0}}{n}\right)D_{\alpha}\\&-3\frac{U_{\gamma\gamma}^{(k)}f_{k0}}{n}G_{\alpha}-\frac{np-2m+2}{4\left(n+m-1\right)}f_{\alpha}\mid\partial f\mid^{4}]\}-e^{2(n-1)f}B\}\\&=e^{\left(\frac{n}{n+m-1}(2+p)+2n-2\right)f}\left(\mid E_{\alpha\overline{\beta}}\mid^{2}+\mid D_{\alpha\beta}\mid^{2}\right)\\&+e^{2(n-1)f}\left(\mid G_{\alpha}\mid^{2}+\mid D_{\alpha\beta}f_{\overline{\gamma}}+E_{\alpha\overline{\gamma}}f_{\beta}\mid^{2}\right)\\&+\left(\mid G_\alpha+D_\alpha\mid^2+\mid G_\alpha-E_\alpha\mid^2\right)e^{2(n-1)f}\\&+e^{(2n-2)f}\mathbf{Re}\left(E_\alpha+D_\alpha\right)f_{\overline{\alpha}}\left(\frac{np-2m+2}{n+m-1}e^{\frac{n}{n+m-1}(2+p)f}-\frac{np-2m+2}{2\left(n+m-1\right)}|\partial f|^2\right)\\&-\frac{np-2m+2}{n+m-1}(2n-1)\mid\partial f\mid^2e^{\left(\frac{n}{n+m-1}(2+p)+2n-2\right)f}\\&+\left(\frac{np-2m+2}{4\left(n+m-1\right)}\left(n+2\right)+\left(1-2n\right)\frac{np-2m+2}{n+m-1}\right)\mid\partial f\mid^4e^{\left(\frac{n}{n+m-1}\left(2+p\right)+2n-2\right)f}\\&-\frac{np-2m+2}{4\left(n+m-1\right)}n\mid\partial f\mid^6e^{(2n-2)f}-3n\frac{np-2m+2}{n+m-1}\mid f_{k0}\mid^2e^{\left(\frac{n}{n+m-1}(2+p)+2n-2\right)f}.
\end{aligned}
\end{split}
\end{equation}
\end{proposition}
where
\begin{equation}\label{2.9}
\begin{split}
\begin{aligned}
&B=\left(\mid\partial f\mid^{2}+e^{\frac{n}{n+m-1}(2+p)f}\right)\left[2U_{\beta\alpha}^{(k)}f_{k0,\alpha}f_{\overline{\beta}}+2U_{\beta\alpha}^{(k)}f_{k0}f_{\alpha\overline{\beta}}-2\frac{U_{\alpha\alpha}^{(k)}}{n}\left(f_{k0}f_{\beta\overline{\beta}}+f_{k0,\beta}f_{\overline{\beta}}\right)\right]\\
&\quad +2\frac{U_{\gamma\gamma}^{(k)}f_{k0}}{n}\left(-2U_{\beta\alpha}^{(k)}f_{k0,\alpha}f_{\overline{\beta}}+2U_{\beta\alpha}^{(k)}f_{k0}f_{\alpha\overline{\beta}}-2\frac{U_{\alpha\alpha}^{(k)}}{n}\left(f_{k0}f_{\beta\overline{\beta}}-f_{k0,\beta}f_{\overline{\beta}}\right)\right).
\end{aligned}
\end{split}
\end{equation}

\begin{remark}
When $m=1$, that is $\mathbb{G}$ is the Heisenberg group, the (\ref{2.7}) is exactly the crucial identity found by Ma and Ou \cite{OM2023} on the Heisenberg group.
\end{remark}
Now we verify this equality (\ref{2.8}) by applying the previously derived equations and commutation rules.
\begin{proof} Denote
$\mathcal{M}=\mathbf{Re}\left(\mathcal{L}_1+\mathcal{L}_2+\mathcal{L}_3+\mathcal{L}_4-e^{2(n-1)f}B\right)$
with
\begin{eqnarray*}
\mathcal{L}_1&=&Z_{\overline{\alpha}}\left\{\left(g+3\frac{U_{\gamma\gamma}^{(k)}f_{k0}}{n}\right)E_{\alpha}e^{2(n-1)f}\right\},\\
\mathcal{L}_{2}&=&Z_{\overline{\alpha}}\left\{\left(g-\frac{U_{\gamma}^{(k)}f_{k0}}{n}\right)D_{\alpha}e^{2(n-1)f}\right\},\\
\mathcal{L}_{3}&=&Z_{\overline{\alpha}}\left\{-3\frac{U_{\gamma\gamma}^{(k)}f_{k0}}{n}G_{\alpha}e^{2(n-1)f}\right\},\\
\mathcal{L}_{4}&=&Z_{\overline{\alpha}}\left\{-\frac{np-2m+2}{4\left(n+m-1\right)}f_{\alpha}\mid\partial f\mid^{4}e^{2(n-1)f}\right\}.
\end{eqnarray*}
First we compute $\mathcal{L}_1$. By (\ref{2.4}) and using the commutative formulae, we can get the following equation.
\begin{eqnarray}\label{2.10}
\begin{split}
E_{\alpha,\overline{\alpha}} &= f_{\alpha\overline{\beta}\overline{\alpha}}f_{\beta}+f_{\alpha\overline{\beta}}f_{\beta\overline{\alpha}}+g_{\overline{\alpha}}f_{\alpha}+gf_{\alpha\overline{\alpha}}\\
&= f_{\alpha\overline{\alpha}\overline{\beta}}f_\beta+f_{\alpha\overline{\beta}}(f_{\overline{\alpha}\beta}+2U_{\beta\alpha}^{(k)}f_{k0})+g_{\overline{\alpha}}f_\alpha+gf_{\alpha\overline{\alpha}}\\
&= -ng_{\overline{\beta}}f_\beta+f_{\alpha\overline{\beta}}f_{\overline{\alpha}\beta}+2U_{\beta\alpha}^{(k)}f_{k0}f_{\alpha\overline{\beta}}+g_{\overline{\alpha}}f_{\alpha}+gf_{\alpha\overline{\alpha}}\\
&= (1-n)f_\alpha g_{\overline{\alpha}}+(E_{\alpha\overline{\beta}}-g\delta_{\alpha\overline{\beta}})(E_{\overline{\alpha}\beta}-\overline{g}\delta_{\overline{\alpha}\beta})-n|g|^2\\
&\quad +2U_{\beta\alpha}^{(k)}f_{k0}f_{\alpha\overline{\beta}}-2\frac{U_{\beta\beta}^{(k)}f_{k0}}{n}f_{\alpha\overline{\alpha}}.
\end{split}
\end{eqnarray}
There for
\begin{equation}\label{2.11}
 e^{-2(n-1)f}\mathcal{L}_{1}=e^{-2(n-1)f}Z_{\overline{\alpha}}\left\{\left(g+3\frac{U_{\gamma\gamma}^{(k)}f_{k0}}{n}\right)E_{\alpha}e^{2(n-1)f}\right\}
 \end{equation}
  \begin{eqnarray*}\label{2.11}
 &=\left(g+3\frac{U_{\gamma\gamma}^{(k)}f_{k0}}{n}\right)E_{\alpha,\overline{\alpha}}\\&\quad+\left(g_{\overline{\alpha}}+3\frac{U_{\gamma\gamma}^{(k)}f_{k0,\alpha}}{n}\right)E_{\alpha}+2(n-1)\left(g+3\frac{U_{\gamma\gamma}^{(k)}f_{k0}}{n}\right)f_{\overline{\alpha}}E_{\alpha}\\
 &=\left(g+3\frac{U_{\gamma\gamma}^{(k)}f_{k0}}{n}\right)\left(|E_{\alpha\overline{\beta}}|^2+(1-n)f_\alpha g_{\overline{\alpha}}+2U_{\beta\alpha}^{(k)}f_{k0}f_{\alpha\overline{\beta}}-2\frac{U_{\beta\beta}^{(k)}f_{k0}}{n}f_{\alpha\alpha}\right)\\&\quad+g_{\overline{\alpha}}E_\alpha+3\left(-G_{\overline{\alpha}}+\overline{g}f_{\overline{\alpha}}\right)E_\alpha+2(n-1)\left(g+3\frac{U_{\gamma\gamma}^{(k)}f_{k0}}{n}\right)f_{\overline{\alpha}}E_\alpha\\&=\left(g+3\frac{U_{\gamma\gamma}^{(k)}f_{k0}}{n}\right)|E_{\alpha\beta}|^2+\left(g_{\overline{\alpha}}-3G_{\overline{\alpha}}\right)E_{\alpha}\\&\quad+\left(3\overline{g}+2(n-1)(g+3\frac{U_{\gamma\gamma}^{(k)}f_{k0}}{n})\right)f_{\overline{\alpha}}E_\alpha+(1-n)(g+3\frac{U_{\gamma\gamma}^{(k)}f_{k0}}{n})f_\alpha g_{\overline{\alpha}}\\&\quad+(g+3\frac{U_{\gamma\gamma}^{(k)}f_{k0}}{n})\left(2U_{\beta\alpha}^{(k)}f_{k0}f_{\alpha\overline{\beta}}-2\frac{U_{\beta\beta}^{(k)}f_{k0}}{n}f_{\alpha\overline{\alpha}}\right).
\end{eqnarray*}

Next we compute $\mathcal{L}_2$. Also by (\ref{2.4}) and using the commutative formulae, we can get the following equation.
\begin{equation}\label{2.12}
\begin{split}
 \begin{aligned}f_{\alpha\beta\overline{\alpha}}&=f_{\alpha\overline{\alpha}\beta}+2U_{\beta\alpha}^{(k)}f_{k0,\alpha}\\&=(f_{\overline{\alpha}\alpha}+2U_{\alpha\alpha}^{(k)}f_{k0})_{,\beta}+2U_{\beta\alpha}^{(k)}f_{k0,\alpha}\\&=-n\overline{g}_{\beta}+2U_{\alpha\alpha}^{(k)}f_{k0,\beta}+2U_{\beta\alpha}^{(k)}f_{k0,\alpha}.\end{aligned}
\end{split}
\end{equation}
By this and (\ref{2.4}), (\ref{2.5}) we deduce
\begin{equation}\label{2.13}
\begin{split}
 \begin{aligned}\\&D_{\alpha,\overline{\alpha}}\\&=f_{\alpha\beta\overline{\alpha}}f_{\overline{\beta}}+f_{\alpha\beta}f_{\overline{\beta}\overline{\alpha}}-2(\mid\partial f\mid^{2})_{,\overline{\alpha}}f_{\alpha}-2\mid\partial f\mid^{2}f_{\alpha\overline{\alpha}}\\&=\left(-n\overline{g}_{\beta}+2U_{\alpha\alpha}^{(k)}f_{k0,\beta}+2U_{\beta\alpha}^{(k)}f_{k0,\alpha}\right)f_{\overline{\beta}}\\&\quad+(D_{\alpha\beta}+2f_{\alpha}f_{\beta})(D_{\overline{\alpha}\overline{\beta}}+2f_{\overline{\alpha}}f_{\overline{\beta}})\\&\quad-2\left(D_{\overline{\alpha}}+E_{\overline{\alpha}}+\overline{g}f_{\overline{\alpha}}-2f_{\overline{\alpha}}e^{\frac{n}{n+m-1}(2+p)f}\right)f_{\alpha}+2n|\partial f|^2g\\&=\mid D_{\alpha\beta}\mid^2+2f_{\overline{\alpha}}D_\alpha-2f_\alpha E_{\overline{\alpha}}-nf_{\overline{\alpha}}\overline{g}_\alpha+\left(2U_{\alpha\alpha}^{(k)}f_{k0,\beta}+2U_{\beta\alpha}^{(k)}f_{k0,\alpha}\right)f_{\overline{\beta}}\\&\quad+2\left|\partial f\right|^2\left(ng+g\right).\end{aligned}
\end{split}
\end{equation}
There for
\begin{equation}\label{2.14}
\begin{split}
 \begin{aligned}&e^{-2(n-1)f}\mathcal{L}_{2}\\&=e^{-2(n-1)f}Z_{\overline{\alpha}}\left\{\left(g-\frac{U_{\gamma\gamma}^{(k)}f_{k0}}{n}\right)D_{\alpha}e^{2(n-1)f}\right\}\\&=\left(g-\frac{U_{\gamma\gamma}^{(k)}f_{k0}}{n}\right)D_{\alpha,\overline{\alpha}}+\left(g_{\overline{\alpha}}-\frac{U_{\gamma\gamma}^{(k)}f_{k0,\overline{\alpha}}}{n}\right)D_\alpha\\&\quad+2(n-1)\left(g-\frac{U_{\gamma\gamma}^{(k)}f_{k0}}{n}\right)f_{\overline{\alpha}}D_\alpha\\&=\left(g-\frac{U_{\gamma\gamma}^{(k)}f_{k0}}{n}\right)\bigg[ |D_{\alpha\beta}|^2 + 2f_{\overline{\alpha}}D_\alpha - 2f_\alpha E_{\overline{\alpha}} - nf_{\overline{\alpha}}\overline{g}_\alpha \\
&\quad + 2U_{\alpha\alpha}^{(k)}f_{k0,\beta}f_{\overline{\beta}} + 2U_{\beta\alpha}^{(k)}f_{k0,\alpha}f_{\overline{\beta}} + 2|\partial f|^2\left(ng+g\right) \bigg]\\&\quad+\left(g_{\overline{\alpha}}+G_{\overline{\alpha}}-\overline{g}f_{\overline{\alpha}}\right)D_{\alpha}+2(n-1)\left(g-\frac{U_{\gamma\gamma}^{(k)}f_{k0}}{n}\right)f_{\overline{\alpha}}D_{\alpha}\end{aligned}
\end{split}
\end{equation}

\begin{equation*}\label{2.14}
\begin{split}
\begin{aligned}
&=\left(g-\frac{U_{\gamma\gamma}^{(k)}f_{k0}}{n}\right)|D_{\alpha\beta}|^2+(g_{\overline{\alpha}}+G_{\overline{\alpha}})D_\alpha+(2ng-\overline{g}-2U_{\alpha\alpha}^{(k)}f_{k0})f_{\overline{\alpha}}D_\alpha\\
&\quad-2(g-\frac{U_{\gamma\gamma}^{(k)}f_{k0}}{n})f_\alpha E_{\overline{\alpha}}-n(g-\frac{U_{\gamma\gamma}^{(k)}f_{k0}}{n})f_{\overline{\alpha}}\overline{g}_\alpha\\&\quad+(g-\frac{U_{\gamma\gamma}^{(k)}f_{k0}}{n})\left(2U_{\alpha\alpha}^{(k)}f_{k0,\beta}f_{\overline{\beta}}+2U_{\beta\alpha}^{(k)}f_{k0,\alpha}f_{\overline{\beta}}+2|\partial f|^2\left(ng+g\right)\right).\end{aligned}
\end{split}
\end{equation*}

Then we compute $\mathcal{L}_3$. Also by (\ref{2.4}) and using the commutative formulae, we can get the following equation.
\begin{equation}\label{2.15}
\begin{split}
\begin{aligned}G_{\alpha,\overline{\alpha}}&=\frac{U_{\beta\beta}^{(k)}f_{k0,\alpha\overline{\alpha}}}{n}+g_{\overline{\alpha}}f_{\alpha}+gf_{\alpha\overline{\alpha}}\\&=\frac{U_{\beta\beta}^{(k)}f_{\alpha\overline{\alpha},k0}}{n}+g_{\overline{\alpha}}f_{\alpha}+g(f_{\overline{\alpha}\alpha}+2U_{\alpha\alpha}^{(k)}f_{k0})\\&=f_{\alpha}g_{\overline{\alpha}}-U_{\beta\beta}^{(k)}g_{k0}-n\mid g\mid^{2}+2U_{\alpha\alpha}^{(k)}f_{k0}g.\end{aligned}
\end{split}
\end{equation}
There for
\begin{equation}\label{2.16}
\begin{split}
\begin{aligned}&e^{-2(n-1)f}\mathcal{L}_3\\&=e^{-2(n-1)f}Z_{\overline{\alpha}}\left\{-3\frac{U_{\gamma\gamma}^{(k)}f_{k0}}{n}G_\alpha e^{2(n-1)f}\right\}\\&=-3\frac{U_{\gamma\gamma}^{(k)}f_{k0}}{n}G_{\alpha,\overline{\alpha}}-3\frac{U_{\gamma\gamma}^{(k)}f_{k0,\overline{\alpha}}}{n}G_\alpha-6(n-1)\frac{U_{\gamma\gamma}^{(k)}f_{k0}}{n}f_{\overline{\alpha}}G_\alpha\\&=-3\frac{U_{\gamma\gamma}^{(k)}f_{k0}}{n}{\left(f_\alpha g_{\overline{\alpha}}-U_{\beta\beta}^{(k)}g_{k0}-n\mid g\mid^2+2U_{\alpha\alpha}^{(k)}f_{k0}g\right)}\\&+3(G_{\overline{\alpha}}-\overline{g}f_{\overline{\alpha}})G_\alpha-6(n-1)\frac{U_{\gamma\gamma}^{(k)}f_{k0}}{n}f_{\overline{\alpha}}G_\alpha\\&=3\mid G_\alpha\mid^2-3(\overline{g}+2(n-1)\frac{U_{\gamma\gamma}^{(k)}f_{k0}}{n})f_{\overline{\alpha}}G_\alpha-3\frac{U_{\gamma\gamma}^{(k)}f_{k0}}{n}f_\alpha g_{\overline{\alpha}}\\&+3\frac{U_{\gamma\gamma}^{(k)}f_{k0}}{n}U_{\beta\beta}^{(k)}g_{k0}+3n\frac{U_{\gamma\gamma}^{(k)}f_{k0}}{n}|g|^2-6\frac{U_{\gamma\gamma}^{(k)}}{n}U_{\alpha\alpha}^{(k)}|f_{k0}|^2g.\end{aligned}
\end{split}
\end{equation}

Furthermore, we compute $\mathcal{L}_4 $. By (\ref{2.3}) and (\ref{2.5}), we get the following equation.
\begin{equation}\label{2.17}
\begin{split}
  \begin{aligned}&e^{-2(n-1)f}\mathcal{L}_4\\&=e^{-2(n-1)f}Z_{\overline{\alpha}}\left\{-\frac{np-2m+2}{4\left(n+m-1\right)}f_\alpha\mid\partial f\mid^4e^{2(n-1)f}\right\}\\&=-\frac{np-2m+2}{2\left(n+m-1\right)}\left(E_{\overline{\alpha}}+D_{\overline{\alpha}}\right)f_{\alpha}\mid\partial f\mid^{2}-\frac{np-2m+2}{4\left(n+m-1\right)}n\mid\partial f\mid^{6}\\&\quad+\frac{np-2m+2}{4\left(n+m-1\right)}(n+2)\mid\partial f\mid^{4}e^{\frac{n}{n+m-1}(2+p)f}\\&\quad-\frac{np-2m+2}{4\left(n+m-1\right)}(n+2)\frac{U_{\alpha\alpha}^{(k)}f_{k0}}{n}|\partial f|^4.\end{aligned}
\end{split}
\end{equation}

Finally, by  (\ref{2.3}),(\ref{2.6}), (\ref{2.7})(\ref{2.11}), (\ref{2.14}), (\ref{2.16}) and (\ref{2.17}), we obtain
\begin{equation}\label{2.18}
\begin{split}
\begin{aligned}&e^{-2(n-1)f}(\mathcal{L}_1+\mathcal{L}_2+\mathcal{L}_3+\mathcal{L}_4)\\&=\Bigg(g-\frac{U_{\gamma\gamma}^{(k)}f_{k0}}{n}\Bigg)\mid D_{\alpha\beta}\mid^{2}+\Bigg(g+3\frac{U_{\gamma\gamma}^{(k)}f_{k0}}{n}\Bigg)\mid E_{\alpha\beta}\mid^{2}+3\mid G_{\alpha}\mid^{2}\\&\quad+(D_{\overline{\alpha}}+E_{\overline{\alpha}}+2G_{\overline{\alpha}})D_{\alpha}+(D_{\overline{\alpha}}+E_{\overline{\alpha}}-2G_{\overline{\alpha}})E_{\alpha}\\&\quad+\Bigg(\frac{np-2m+2}{n+m-1} f_{\overline{\alpha}}e^{\frac{n}{n+m-1}(2+p)f}-\frac{np-2m+2}{2\left(n+m-1\right)}f_{\overline{\alpha}}\mid\partial f\mid^2\Bigg)\Bigg(D_\alpha+E_\alpha\Bigg)\\&\quad+(n-1)\Bigg(\mid\partial f\mid^{2}+e^{\frac{n}{n+m-1}(2+p)f}\Bigg)\Bigg(f_{\overline{\alpha}}G_{\alpha}-f_{\alpha}G_{\overline{\alpha}}\Bigg)\\&\quad+(5n+1)\frac{U_{\gamma\gamma}^{(k)}f_{k0}}{n}(f_{\overline{\alpha}}G_{\alpha}+f_{\alpha}G_{\overline{\alpha}})\\&\quad-\frac{np-2m+2}{n+m-1}(2n-1)\mid\partial f\mid^2e^{\frac{2n}{n+m-1}(2+p)f}\\
&\quad+\Bigg(\frac{np-2m+2}{4\left(n+m-1\right)}\big(n+2\big)+\big(1-2n\big)\frac{np-2m+2}{n+m-1}\Bigg)\mid\partial f\mid^4 e^{\frac{n}{n+m-1}(2+p)f}\end{aligned}
\end{split}
\end{equation}
\begin{equation*}\label{2.18}
\begin{split}
\begin{aligned}
&\quad-\frac{np-2m+2}{4\left(n+m-1\right)}n\mid\partial f\mid^6-3n\frac{np-2m+2}{n+m-1}\mid f_{k0}\mid^2e^{\frac n{n+m-1}(2+p)f}\\&\quad-\frac{np-2m+2}{4\left(n+m-1\right)}(n+2)\frac{U_{\alpha\alpha}^{(k)}f_{k0}}n\mid\partial f\mid^4\\&\quad-\frac{np-2m+2}{n+m-1}\frac{U_{\alpha\alpha}^{(k)}f_{k0}}n\mid\partial f\mid^2e^{\frac n{n+m-1}(2+p)f}\\&\quad+3n\frac{U_{\gamma\gamma}^{(k)}f_{k0}}{n}| g |^{2}+6\frac{U_{\gamma\gamma}^{(k)}U_{\beta\beta}^{(k)}}{n^{2}}U_{\alpha\alpha}^{(k)}| f_{k0} |^{2} f_{k0}+3n\frac{U_{\gamma\gamma}^{(k)}f_{k0}}{n}f_{k0,k_{1}0}\\&\quad+\left(\mid\partial f\mid^2+e^{\frac{n}{n+m-1}(2+p)f}\right)\left[\left(2U_{\beta\alpha}^{(k)}f_{k0,\alpha}f_{\beta}+2U_{\beta\alpha}^{(k)}f_{k0}f_{\alpha\overline{\beta}}\right)-2\frac{U_{\alpha\alpha}^{(k)}}{n}\left(f_{k0}f_{\beta\overline{\beta}}+f_{k0,\beta}f_{\overline{\beta}}\right)\right]\\&\quad+2\frac{U_{\gamma\gamma}^{(k)}f_{k0}}{n}\left[\left(-2U_{\beta\alpha}^{(k)}f_{k0,\alpha}f_{\overline{\beta}}+2U_{\beta\alpha}^{(k)}f_{k0}f_{\alpha\overline{\beta}}\right)-2\frac{U_{\alpha\alpha}^{(k)}}{n}\left(f_{k0}f_{\beta\overline{\beta}}-f_{k0,\beta}f_{\overline{\beta}}\right)\right].\end{aligned}
\end{split}
\end{equation*}
The last two items are $B$. So, since $f$ is real, we can complete the proof of proposition \ref{Pro-1} easily by some simplification.

Next, we verify some properties of $\mathcal{M}$. By (\ref{2.17}), we obtain
\begin{equation}\label{2.19}
\begin{split}
\begin{aligned}
&\mathcal{M}=e^{\left(\frac{n}{n+m-1}(2+p)+2n-2\right)f}\left(\mid E_{\alpha\overline{\beta}}\mid^{2}+\mid D_{\alpha\beta}\mid^{2}\right)\\&+e^{2(n-1)f}\left(\mid G_{\alpha}\mid^{2}+\mid D_{\alpha\beta}f_{\overline{\gamma}}+E_{\alpha\overline{\gamma}}f_{\beta}\mid^{2}\right)\\&+s\left(\mid G_\alpha+D_\alpha\mid^2+\mid G_\alpha-E_\alpha\mid^2\right)e^{2(n-1)f}\\&+\left(1-s\right)\left(\mid G_\alpha+D_\alpha\mid^2+\mid G_\alpha-E_\alpha\mid^2\right)e^{2(n-1)f}\\&+e^{(2n-2)f}\mathbf{Re}\left(G_\alpha+D_\alpha\right)f_{\overline{\alpha}}\left(\frac{np-2m+2}{n+m-1}e^{\frac{n}{n+m-1}(2+p)f}-\frac{np-2m+2}{2\left(n+m-1\right)}|\partial f|^2\right)\\&+e^{(2n-2)f}\mathbf{Re}\left(E_\alpha-G_\alpha\right)f_{\overline{\alpha}}\left(\frac{np-2m+2}{n+m-1}e^{\frac{n}{n+m-1}(2+p)f}-\frac{np-2m+2}{2\left(n+m-1\right)}|\partial f|^2\right)\\&-\frac{np-2m+2}{n+m-1}(2n-1)\mid\partial f\mid^2e^{\left(\frac{n}{n+m-1}(2+p)+2n-2\right)f}\\&+\left(\frac{np-2m+2}{4\left(n+m-1\right)}\left(n+2\right)+\left(1-2n\right)\frac{np-2m+2}{n+m-1}\right)\mid\partial f\mid^4e^{\left(\frac{n}{n+m-1}\left(2+p\right)+2n-2\right)f}\\&-\frac{np-2m+2}{4\left(n+m-1\right)}n\mid\partial f\mid^6e^{(2n-2)f}-3n\frac{np-2m+2}{n+m-1}\mid f_{k0}\mid^2e^{\left(\frac{n}{n+m-1}(2+p)+2n-2\right)f}.
\end{aligned}
\end{split}
\end{equation}
To verify that every term in $\mathcal{M}$ is nonnegative, we rewrite \label{2.18} as
\begin{equation}\label{2.20}
\begin{split}
\begin{aligned}&\mathcal{M}=e^{\left(\frac{n}{n+m-1}(2+p)+2n-2\right)f}\left(\mid E_{\alpha\overline{\beta}}\mid^2+\mid D_{\alpha\beta}\mid^2\right)+e^{2(n-1)f}\left(\mid G_\alpha\mid^2+\mid D_{\alpha\beta}f_{\overline{\gamma}}+E_{\alpha\overline{\gamma}}f_\beta\mid^2\right)\\&+s\left(\mid G_\alpha+D_\alpha\mid^2+\mid G_\alpha-E_\alpha\mid^2\right)e^{2(n-1)f}\\&+e^{(2n-2)f}\left|\sqrt{1-s}\left(G_\alpha+D_\alpha\right)+\frac{1}{2\sqrt{1-s}}f_\alpha\left(\frac{np-2m+2}{n+m-1}e^{\frac{n}{n+m-1}(2+p)f}-\frac{np-2m+2}{2(n+m-1)}|\partial f|^2\right)\right|^2\\&+e^{(2n-2)f}\left|\sqrt{1-s}\left(E_\alpha-G_\alpha\right)+\frac{1}{2\sqrt{1-s}}f_\alpha\left(\frac{np-2m+2}{n+m-1}e^{\frac{n}{n+m-1}(2+p)f}-\frac{np-2m+2}{2\left(n+m-1\right)}|\partial f|^2\right)\right|^2\\
&-\frac{1}{2\left(1-s\right)}e^{(2n-2)f}\left|\partial f\right|^2\left[\left(\frac{np-2m+2}{n+m-1}\right)^2e^{\frac{n}{n+m-1}2(2+p)f}+\frac{1}{4}\left(\frac{np-2m+2}{n+m-1}\right)^2\left|\partial f\right|^4\right]\\&-\left(\frac{np-2m+2}{n+m-1}\right)^2e^{\frac{n}{n+m-1}(2+p)f}\left|\partial f\right|^2\\&-\frac{np-2m+2}{n+m-1}(2n-1)\mid\partial f\mid^2e^{\left(\frac{2n}{n+m-1}(2+p)+2n-2\right)f}\\&+\left(\frac{np-2m+2}{4\left(n+m-1\right)}(n+2)+(1-2n)\frac{np-2m+2}{n+m-1}\right)\mid\partial f\mid^4e^{\left(\frac{n}{n+m-1}(2+p)+2n-2\right)f}\\&-\frac{np-2m+2}{4(n+m-1)}n\mid\partial f\mid^6e^{(2n-2)f}-3n\frac{np-2m+2}{n+m-1}\mid f_{k0}\mid^2e^{\left(\frac{n}{n+m-1}(2+p)+2n-2\right)f}.\end{aligned}
\end{split}
\end{equation}
Continue to get
\begin{equation}\label{2.21}
\begin{split}
\begin{aligned}&\mathcal{M}=e^{\left(\frac{n}{n+m-1}(2+p)+2n-2\right){f}}\left(\mid E_{\alpha\overline{\beta}}\mid^{2}+\mid D_{\alpha\beta}\mid^{2}\right)+e^{2(n-1)f}\left(\mid G_{\alpha}\mid^{2}+\mid D_{\alpha\beta}f_{\overline{\gamma}}+E_{\alpha\overline{\gamma}}f_{\beta}\mid^{2}\right)\\&+s\left(\mid G_\alpha+D_\alpha\mid^2+\mid G_\alpha-E_\alpha\mid^2\right)e^{2(n-1)f}\\&+e^{(2n-2)f}\left|\sqrt{1-s}\left(G_\alpha+D_\alpha\right)+\frac{1}{2\sqrt{1-s}}f_\alpha\left(\frac{np-2m+2}{n+m-1}e^{\frac{n}{n+m-1}(2+p)f}-\frac{np-2m+2}{2(n+m-1)}|\partial f|^2\right)\right|^2\\
&+e^{(2n-2)f}\left|\sqrt{1-s}\left(E_\alpha-G_\alpha\right)+\frac{1}{2\sqrt{1-s}}f_\alpha\left(\frac{np-2m+2}{n+m-1}e^{\frac{n}{n+m-1}(2+p)f}-\frac{np-2m+2}{2(n+m-1)}|\partial f|^2\right)\right|^2\end{aligned}
\end{split}
\end{equation}
\begin{equation*}\label{2.21}
\begin{split}
\begin{aligned}
&-\left(\frac{np-2m+2}{n+m-1}(2n-1)+\frac{1}{2\left(1-s\right)}\left(\frac{np-2m+2}{n+m-1}\right)^2\right)|\partial f|^2e^{\left(\frac{2n}{n+m-1}(2+p)+2n-2\right)f}\\&+\bigg[\frac{np-2m+2}{4\left(n+m-1\right)}(n+2)+(1-2n)\frac{np-2m+2}{n+m-1}\\&+\frac{1}{2\left(1-s\right)}\left(\frac{np-2m+2}{n+m-1}\right)^2\bigg]\mid\partial f\mid^4e^{\left(\frac{n}{n+m-1}(2+p)+2n-2\right)}\\&-\left(\frac{np-2m+2}{4\left(n+m-1\right)}n+\frac{1}{4}\left(\frac{np-2m+2}{n+m-1}\right)^2\frac{1}{2\left(1-s\right)}\right)\mid\partial f\mid^6e^{(2n-2)f}\\&-3n\frac{np-2m+2}{n+m-1}|f_{k0}|^2e^{\left(\frac{n}{n+m-1}(2+p)+2n-2\right)f}.\end{aligned}
\end{split}
\end{equation*}
Now we take $0<s=s_0=\frac{1}{2}+\frac{np-2m+2}{4n\left(n+m-1\right)}<1$, then
\begin{equation}\label{2.22}
\begin{split}
\begin{aligned}&\mathcal{M}=e^{\left(\frac{n}{n+m-1}(2+p)+2n-2\right)f}\left(\mid E_{\alpha\overline{\beta}}\mid^2+\mid D_{\alpha\beta}\mid^2\right)+e^{2(n-1)f}\left(\mid G_\alpha\mid^2+\mid D_{\alpha\beta}f_{\overline{\gamma}}+E_{\alpha\overline{\gamma}}f_\beta\mid^2\right)\\&+s\left(\mid G_\alpha+D_\alpha\mid^2+\mid G_\alpha-E_\alpha\mid^2\right)e^{2(n-1)f}\\&+e^{(2n-2)f}\left|\sqrt{1-s}\left(G_\alpha+D_\alpha\right)+\frac{1}{2\sqrt{1-s}}f_\alpha\left(\frac{np-2m+2}{n+m-1}e^{\frac{n}{n+m-1}(2+p)f}-\frac{np-2m+2}{2(n+m-1)}|\partial f|^2\right)\right|^2\\&+e^{(2n-2)f}\left|\sqrt[]{1-s}\left(E_\alpha-G_\alpha\right)+\frac{1}{2\sqrt{1-s}}f_\alpha\left(\frac{np-2m+2}{n+m-1}e^{\frac{n}{n+m-1}(2+p)f}-\frac{np-2m+2}{2(n+m-1)}|\partial f|^2\right)\right|^2\\&-\frac{np-2m+2}{n+m-1}\frac{\left(2n^2-n\right)\left(2n+2m-2\right)+\left(np-2m+2\right)}{2n\left(n+m-1\right)-\left(np-2m+2\right)}\mid\partial f\mid^2e^{\left(\frac{2n}{n+m-1}\left(2+p\right)+2n-2\right)f}\\&+\frac{np-2m+2}{n+m-1}\frac{\left(\frac{3}{2}-\frac{7}{4}n\right)2n\left(n+m-1\right)+\left(\frac{1}{4}n-\frac{3}{2}\right)(np-2m+2)}{2n\left(n+m-1\right)-\left(mp-2m+2\right)}|\partial f|^4e^{\left(\frac{n}{n+m-1}(2+p)+2n-2\right)f}\\&-\frac{np-2m+2}{4\left(n+m-1\right)}\frac{n^2\left(2n+2m-2\right)+n\left(np-2m+2\right)}{2n\left(n+m-1\right)-\left(np-2m+2\right)}\mid\partial f\mid^6e^{(2n-2)f}\\&-3n\frac{np-2m+2}{n+m-1}|f_{k0}|^2e^{\left(\frac{n}{n+m-1}(2+p)+2n-2\right)f}.\end{aligned}
\end{split}
\end{equation}
Now clearly, all the coefficients above are positive for $-2<p<0$ and  $\mathcal{M}\geq 0$.\end{proof}

\section{Proof of Theorem \ref{Thm2} }

\setcounter{equation}{0}
\setcounter{theorem}{0}

Now we will give the proof of Theorem \ref{Thm2} based on the conclusion in Section 2.
Since $B_{4r}\subset \Omega$, we can take  a real smooth cut off function  $\eta$  such that

\begin{equation}\label{3-1}
\begin{cases}
                          \eta\equiv 1  &in \,\,B_r,\\
                        0\leq\eta\leq1  &in \,\,B_{2r},\\
                          \eta\equiv 0  &in \,\,\Omega\backslash B_{2r},\\
     |\partial \eta|\lesssim \frac{1}{r}  &in \,\,\Omega.
\end{cases}
\end{equation}
\noindent In order to make the proof process more concise where here and in the sequel we use ``$\lesssim $" , ``$\cong$"  to replace ``$\leq$" and ``$=$" respectively, to eliminate the effect of some positive constants independent of $r$ and $f$.

Rewrite (\ref{2.8}) and take a real $s>0$ big enough. Multiply both sides of it by $\eta^s $ and integrate over  $\Omega$ we have
\begin{equation}\label{3-2}
\begin{split}
\begin{aligned}
&\int_{\Omega}\eta^{s}\mathcal{M}\\
&=\int_{\Omega}\eta^{s}\mathbf{Re}\Bigg\{Z_{\overline{\alpha}}\Bigg\{\left[\left(D_{\alpha}+E_{\alpha}\right)\left(\mid\partial f\mid^{2}+e^{\frac{m}{n+m-1}(2+p)f}\right)\right.\\
&\quad\left.-\frac{U_{\gamma\gamma}^{(k)}f_{k0}}{n}(2D_{\alpha}-2E_{\alpha}+3G_{\alpha})-\frac{np-2m+2}{4\left(n+m-1\right)}f_{\alpha}\mid\partial f\mid^{4}\right]e^{2(n-1)f}\Bigg\}\\
&\quad-e^{2(n-1)f}B\Bigg\}.
\end{aligned}
\end{split}
\end{equation}
Integrating by part and using (\ref{3-1}) we get
\begin{equation}\label{3-3}
\begin{split}
\begin{aligned}&\int_\Omega\eta^s\mathcal{M}\\&=-s\int_\Omega\eta^{s-1}\mathbf{Re}\eta_{\overline{\alpha}}\{[(D_\alpha+E_\alpha)(\mid\partial f\mid^2+e^{\frac{n}{n+m-1}(2+p)f})\\&-\frac{U_{\gamma\gamma}^{(k)}f_{k0}}{n}(2D_{\alpha}-2E_{\alpha}+3G_{\alpha})-\frac{np-2m+2}{4\left(n+m-1\right)}f_{\alpha}\mid\partial f\mid^{4}]e^{2(n-1)f}\}\\&-\int_\Omega\eta^s\mathbf{Re}e^{2(n-1)f}B\\&\lesssim\frac{1}{r}\int_{\Omega}\eta^{s-1}\left\{\mid D_{\alpha}+E_{\alpha}\mid(\mid\partial f\mid^{2}+e^{\frac{n}{n+m-1}(2+p)f})e^{2(n-1)f}\right.\\&+\mid f_{k0}\parallel2D_{\alpha}-2E_{\alpha}+3G_{\alpha}\mid e^{2(n-1)f}+\mid\partial f\mid^{5}e^{2(n-1)f}\}-\int_{\Omega}\eta^{s}\mathbf{Re}e^{2(n-1)f}B.\end{aligned}
\end{split}
\end{equation}
Since
$$|D_{\alpha}+E_{\alpha}|\leq |D_{\alpha}+G_{\alpha}|+|E_{\alpha}-G_{\alpha}|,$$

$$\big|2D_{\alpha}-2E_{\alpha}+3 G_{\alpha}\big| \leq 2|D_{\alpha}+G_{\alpha}|+2|E_{\alpha}-G_{\alpha}|+|G_{\alpha}|,$$
\noindent combine this and Young's inequality for $p>1,q>1,\frac{1}{p}+\frac{1}{q}=1,$ then $ab\leq \epsilon a^p+\frac{(p\varepsilon)^{-\frac{q}{p}}}{q}b^q$ in (\ref{3-3}) we obtain
\begin{equation}\label{3-4}
\begin{split}
\begin{aligned}&\int_{\Omega}\eta^{s}\mathcal{M}\\&\lesssim\epsilon\int_{\Omega}\eta^{s-1}\left(\mid D_{\alpha}+G_{\alpha}\mid^{2}+\mid E_{\alpha}-G_{\alpha}\mid^{2}+\mid G_{\alpha}\mid^{2}\right)e^{2(n-2)f}\\&+\frac{1}{r^2}\int_\Omega\eta^{s-1}\left(\mid\partial f\mid^4+e^{\frac{2n}{n+m-1}(2+p)f}+\mid f_0\mid^2\right)e^{2nf}\\&+\frac{1}{r}\int_\Omega\eta^{s-1}|\partial f|^5e^{2(n-1)f}-\int_\Omega\eta^s\mathbf{Re}e^{2(n-1)f}B,\end{aligned}
\end{split}
\end{equation}
taking $\epsilon$ small, we have
\begin{equation}\label{3-5}\begin{aligned}&\int_{\Omega}\eta^{s}\mathcal{M}\\&\lesssim\frac{1}{r^2}\int_\Omega\eta^{s-1}\left(\mid\partial f\mid^4+e^{\frac{2n}{n+m-1}(2+p)f}+\mid f_0\mid^2\right)e^{2nf}\\&+\frac{1}{r}\int_\Omega\eta^{s-1}|\partial f|^5e^{2(n-1)f}-\int_\Omega\eta^s\mathbf{Re}e^{2(n-1)f}B,\end{aligned}
\end{equation}
then for the first two terms continue to use Young's inequality we have
\begin{equation}\label{3-6}\begin{aligned}&\frac{1}{r^2}\int_\Omega\eta^{s-1}\mid\partial f\mid^4e^{2nf}\\&\lesssim\epsilon\int_\Omega\eta^{s-1}\mid\partial f\mid^6e^{2nf}+\frac{1}{r^6}\int_\Omega\eta^{s-1}e^{2nf},\end{aligned}
\end{equation}
\begin{equation}\label{3-7}\begin{aligned}&\frac{1}{r^2}\int_\Omega\eta^{s-1}\mid f_0\mid^2e^{2nf}\\&\lesssim\epsilon\int_\Omega\eta^{s-1}\mid f_0\mid^3e^{2nf}+\frac{1}{r^6}\int_\Omega\eta^{s-1}e^{2nf},\end{aligned}
\end{equation}
and
\begin{equation}\label{3-8}\begin{aligned}&\frac{1}{r}\int_\Omega\eta^{s-1}\mid \partial{f}\mid^5e^{2nf}\\&\lesssim\epsilon\int_\Omega\eta^{s-1}|\partial f|^6e^{(2n-\frac{12}{5})f}+\frac{1}{r^6}\int_\Omega\eta^{s-1}e^{2nf}.\end{aligned}
\end{equation}
Inserting these into (\ref{3-5}) and taking $\epsilon$ small we get
\begin{equation}\label{3-9}
\begin{split}
\begin{aligned}&\int_{\Omega}\eta^{s}\mathcal{M}\\&\lesssim\frac{1}{r^2}\int_\Omega\eta^{s-1}e^{2nf}e^{\frac{2n}{n+m-1}(2+p)f}\\&+\frac{1}{r^6}\int_\Omega\eta^{s-1}e^{2nf}-\int_\Omega\eta^s\mathbf{Re}e^{2(n-1)f}B.\end{aligned}
\end{split}
\end{equation}

Now we deal with the last term $-\int_\Omega\eta^s\mathbf{Re}e^{2(n-1)f}B$
\begin{equation}\label{3-10}\begin{aligned}&-\int_\Omega\eta^s\mathbf{Re}e^{2(n-1)f}B\\&=\int_\Omega s\eta^{s-1}\mathbf{Re}\eta_\alpha e^{2(n-1)f}\left[\left(\mid\partial f\mid^2+e^{\frac{n}{n+m-1}(2+p)f}-2\frac{U_{\gamma\gamma}^{(k)}f_{k0}}{n}\right)\left(2U_{\beta\alpha}^{(k)}f_{\overline{\beta}}\right)\right] f_{k0}\\&+\int_\Omega\eta^s\mathbf{Re}\left[e^{2(n-1)f}\left(\mid\partial f\mid^2+e^{\frac{n}{n+m-1}(2+p)f}-2\frac{U_{\gamma\gamma}^{(k)}f_{k0}}{n}\right)\left(2U_{\beta\alpha}^{(k)}f_{\overline{\beta}}\right)\right]_\alpha f_{k0}\\&+\int_\Omega s\eta^{s-1}\mathbf{Re}\eta_{\overline{\beta}}e^{2(n-1)f}\left[\left(\mid\partial f\mid^2+e^{\frac{n}{n+m-1}(2+p)f}+2\frac{U_{\gamma\gamma}^{(k)}f_{k0}}{n}\right)\left(2U_{\beta\alpha}^{(k)}f_{k0}\right)\right]f_\alpha\\&+\int_\Omega\eta^s\mathbf{Re}\left[e^{2(n-1)f}\left(\mid\partial f\mid^2+e^{\frac{n}{n+m-1}(2+p)f}+2\frac{U_{\gamma\gamma}^{(k)}f_{k0}}{n}\right)\left(2U_{\beta\alpha}^{(k)}f_{k0}\right)\right]_{\overline{\beta}}f_\alpha\\&+\int_\Omega s\eta^{s-1}\mathbf{Re}\eta_{\overline{\beta}}e^{2(n-1)f}\left[\left(\mid\partial f\mid^2+e^{\frac{n}{n+m-1}(2+p)f}+2\frac{U_{\gamma\gamma}^{(k)}f_{k0}}{n}\right)\left(-2\frac{U_{\alpha\alpha}^{(k)}}{n}f_{k0}\right)\right]f_\beta\\&+\int_\Omega\eta^s\mathbf{Re}\left[e^{2(n-1)f}\left(\mid\partial f\mid^2+e^{\frac{n}{n+m-1}(2+p)f}+2\frac{U_{\gamma\gamma}^{(k)}f_{k0}}{n}\right)\left(-2\frac{U_{\alpha\alpha}^{(k)}}{n}f_{k0}\right)\right]_{\overline{\beta}}f_\beta\\&+\int_\Omega s\eta^{s-1}\mathbf{Re}\eta_\beta e^{2(n-1)f}\left[\left(\mid\partial f\mid^2+e^{\frac{n}{n+m-1}(2+p)f}-2\frac{U_{\gamma\gamma}^{(k)}f_{k0}}{n}\right)\left(-2\frac{U_{\alpha\alpha}^{(k)}}{n}f_{\overline{\beta}}\right)\right]f_{k0}\\&+\int_\Omega\eta^s\mathbf{Re}\left[e^{2(n-1)f}\left(\mid\partial f\mid^2+e^{\frac{n}{n+m-1}(2+p)f}-2\frac{U_{\gamma\gamma}^{(k)}f_{k0}}{n}\right)\left(-2\frac{U_{\alpha\alpha}^{(k)}}{n}f_{k0}\right)\right]_\beta f_{k0}.\end{aligned}\end{equation}
Use Young's inequality and partial integration to handle each of the above terms, for the first line
\begin{equation}\label{3-11}\begin{aligned}&\int_\Omega s\eta^{s-1}\mathbf{Re}\eta_\alpha e^{2(n-1)f}\left[\left(\mid\partial f\mid^2+e^{\frac{n}{n+m-1}(2+p)f}-2\frac{U_{\gamma\gamma}^{(k)}f_{k0}}{n}\right)\left(2U_{\beta\alpha}^{(k)}f_{\overline{\beta}}\right)\right]f_{k0}\\&\leq\int_\Omega s\eta^{s-1}\mid\eta_\alpha\mid e^{2(n-1)f}\left(\mid\partial f\mid^2+e^{\frac{n}{n+m-1}(2+p)f}+\left|2\frac{U_{\gamma\gamma}^{(k)}f_{k0}}{n}\right|\right)\mid{2U_{\beta\alpha}^{(k)}f_{\overline{\beta}}}\mid\mid{f_{k0}}\mid\\&\lesssim\frac{1}{r}\int_\Omega\eta^{s-1}e^{2(n-1)f}e^{\frac{n}{n+m-1}(2+p)f}\mid{f_{\overline{\beta}}}\mid\mid{f_{k0}}\mid\\&\quad+\frac{1}{r}\int_\Omega \eta^{s-1}e^{2(n-1)f}\mid{f_{k0}}\mid^2\mid{f_{\overline{\beta}}}\mid\\&\quad+\frac{1}{r}\int_\Omega \eta^{s-1}e^{2(n-1)f}\mid\partial f\mid^2\mid{f_{\overline{\beta}}}\mid\mid{f_{k0}}\mid,\end{aligned}\end{equation}
for above three terms using Young's inequality we get
\begin{equation}\label{3-12}\begin{aligned}&\frac{1}{r}\int_\Omega\eta^{s-1}e^{2(n-1)f}e^{\frac{n}{n+m-1}(2+p)f}\mid{f_{\overline{\beta}}}\mid\mid{f_{k0}}\mid\\&\lesssim\epsilon\int_\Omega\eta^{s-1}e^{(2n-\frac{8}{3})f}e^{\frac{n}{n+m-1}(2+p)f}\mid{f_{\overline{\beta}}}\mid^\frac{4}{3}\mid{f_{k0}}\mid^\frac{4}{3}\\&\quad+\frac{1}{r^4}\int_\Omega\eta^{s-1}e^{2nf}e^{\frac{n}{n+m-1}(2+p)f},\end{aligned}\end{equation}

\begin{equation}\label{3-13}\begin{aligned}&\frac{1}{r}\int_\Omega \eta^{s-1}e^{2(n-1)f}\mid{f_{k0}}\mid^2\mid{f_{\overline{\beta}}}\mid\\&\lesssim\epsilon\int_\Omega\eta^{s-1}e^{(2n-\frac{12}{5})f}\mid{f_{\overline{\beta}}}\mid^\frac{6}{5}\mid{f_{k0}}\mid^\frac{12}{5}+\frac{1}{r^6}\int_\Omega\eta^{s-1}e^{2nf}\end{aligned}\end{equation}
and
\begin{equation}\label{3-14}\begin{aligned}&\frac{1}{r}\int_\Omega \eta^{s-1}e^{2(n-1)f}\mid\partial f\mid^2\mid{f_{\overline{\beta}}}\mid\mid{f_{k0}}\mid\\&\lesssim\epsilon\int_\Omega\eta^{s-1}e^{(2n-\frac{12}{5})f}\mid{f_{\overline{\beta}}}\mid^\frac{6}{5}\mid{f_{k0}}\mid^\frac{6}{5}\mid\partial f\mid^\frac{12}{5}+\frac{1}{r^6}\int_\Omega\eta^{s-1}e^{2nf}.\end{aligned}\end{equation}
Inserting these into (\ref{3-11}) and taking $\epsilon$ small we get
\begin{equation}\label{3-15}\begin{aligned}&\int_{\Omega}s\eta^{s-1}\mathbf{Re}\eta_{\alpha}e^{2(n-1)f}\left[\left(\mid\partial f\mid^{2}+e^{\frac{n}{n+m-1}(2+p)f}-2\frac{U_{\gamma\gamma}^{(k)}f_{k0}}{n}\right)\left(2U_{\beta\alpha}^{(k)}f_{\overline{\beta}}\right)\right] f_{k0}\\&\lesssim\frac{1}{r^{4}}{\int}_{\Omega}\eta^{s-1}e^{2nf}e^{\frac{n}{n+m-1}(2+p)f}+\frac{1}{r^{6}}{\int}_{\Omega}\eta^{s-1}e^{2nf}.\end{aligned}\end{equation}
For the second line, using Young's inequality we get
\begin{equation}\label{3-16}\begin{aligned}&\int_\Omega\eta^s\mathbf{Re}\left[ e^{2(n-1)f}\left(\mid\partial f\mid^2+e^{\frac{n}{n+m-1}(2+p)f}-2\frac{U_{\gamma\gamma}^{(k)}f_{k0}}{n}\right)\left(2U_{\beta\alpha}^{(k)}f_{\overline{\beta}}\right)\right]_\alpha f_{k0}\\&\leq\int_{\Omega}\eta^s\mid f_{k0}\mid\mathbf{Re}\left[e^{2(n-1)f}\left(\mid\partial f\mid^2+e^{\frac{n}{n+m-1}(2+p)f}-2\frac{U_{\gamma\gamma}^{(k)}f_{k0}}{n}\right)\left(2U_{\beta\alpha}^{(k)}f_{\overline{\beta}}\right)\right]_\alpha\\&\lesssim\epsilon\int_{\Omega}\eta^s\mid f_{k0}\mid^2\mathbf{Re}\left[e^{2(n-1)f}\left(\mid\partial f\mid^2+e^{\frac{n}{n+m-1}(2+p)f}-2\frac{U_{\gamma\gamma}^{(k)}f_{k0}}{n}\right)\left(2U_{\beta\alpha}^{(k)}f_{\overline{\beta}}\right)\right]_\alpha\\&\quad+\int_{\Omega}\eta^s\mathbf{Re}\left[e^{2(n-1)f}\left(\mid\partial f\mid^2+e^{\frac{n}{n+m-1}(2+p)f}-2\frac{U_{\gamma\gamma}^{(k)}f_{k0}}{n}\right)\left(2U_{\beta\alpha}^{(k)}f_{\overline{\beta}}\right)\right]_\alpha.\end{aligned}\end{equation}
Taking $\epsilon$ small and integrating by part we get
\begin{equation}\label{3-17}\begin{aligned}&\int_\Omega\eta^s\mathbf{Re}\left[ e^{2(n-1)f}\left(\mid\partial f\mid^2+e^{\frac{n}{n+m-1}(2+p)f}-2\frac{U_{\gamma\gamma}^{(k)}f_{k0}}{n}\right)\left(2U_{\beta\alpha}^{(k)}f_{\overline{\beta}}\right)\right]_\alpha f_{k0}\\&\lesssim\int_{\Omega}\eta^s\mathbf{Re}\left[e^{2(n-1)f}\left(\mid\partial f\mid^2+e^{\frac{n}{n+m-1}(2+p)f}-2\frac{U_{\gamma\gamma}^{(k)}f_{k0}}{n}\right)\left(2U_{\beta\alpha}^{(k)}f_{\overline{\beta}}\right)\right]_\alpha\\&=\int_{\Omega}\eta^{s-1}\mathbf{Re}\eta_\alpha\left[e^{2(n-1)f}\left(\mid\partial f\mid^2+e^{\frac{n}{n+m-1}(2+p)f}-2\frac{U_{\gamma\gamma}^{(k)}f_{k0}}{n}\right)\left(2U_{\beta\alpha}^{(k)}f_{\overline{\beta}}\right)\right]\\&\lesssim\frac{1}{r}\int_{\Omega}\eta^{s-1}\mathbf{Re}\left[e^{2(n-1)f}\left(\mid\partial f\mid^2+e^{\frac{n}{n+m-1}(2+p)f}-2\frac{U_{\gamma\gamma}^{(k)}f_{k0}}{n}\right)\left(2U_{\beta\alpha}^{(k)}f_{\overline{\beta}}\right)\right]\\&\lesssim\frac{1}{r}\int_\Omega\eta^{s-1}e^{2(n-1)f}e^{\frac{n}{n+m-1}(2+p)f}\mid{f_{\overline{\beta}}}\mid\\&\quad+\frac{1}{r}\int_\Omega \eta^{s-1}e^{2(n-1)f}\mid{f_{k0}}\mid\mid{f_{\overline{\beta}}}\mid\\&\quad+\frac{1}{r}\int_\Omega \eta^{s-1}e^{2(n-1)f}\mid\partial f\mid^2\mid{f_{\overline{\beta}}}\mid.\end{aligned}\end{equation}
for above three terms using Young's inequality we get
\begin{equation}\label{3-18}\begin{aligned}&\frac{1}{r}\int_\Omega\eta^{s-1}e^{2(n-1)f}e^{\frac{n}{n+m-1}(2+p)f}\mid{f_{\overline{\beta}}}\mid\\&\lesssim\epsilon\int_\Omega\eta^{s-1}e^{(2n-\frac{8}{3})f}e^{\frac{n}{n+m-1}(2+p)f}\mid{f_{\overline{\beta}}}\mid^\frac{4}{3}\\&\quad+\frac{1}{r^4}\int_\Omega\eta^{s-1}e^{2nf}e^{\frac{n}{n+m-1}(2+p)f},\end{aligned}\end{equation}

\begin{equation}\label{3-19}\begin{aligned}&\frac{1}{r}\int_\Omega \eta^{s-1}e^{2(n-1)f}\mid{f_{k0}}\mid\mid{f_{\overline{\beta}}}\mid\\&\lesssim\epsilon\int_\Omega\eta^{s-1}e^{(2n-\frac{12}{5})f}\mid{f_{\overline{\beta}}}\mid^\frac{6}{5}\mid{f_{k0}}\mid^\frac{6}{5}+\frac{1}{r^6}\int_\Omega\eta^{s-1}e^{2nf}\end{aligned}\end{equation}
and
\begin{equation}\label{3-20}\begin{aligned}&\frac{1}{r}\int_\Omega \eta^{s-1}e^{2(n-1)f}\mid\partial f\mid^2\mid{f_{\overline{\beta}}}\mid\\&\lesssim\epsilon\int_\Omega\eta^{s-1}e^{(2n-\frac{12}{5})f}\mid{f_{\overline{\beta}}}\mid^\frac{6}{5}\mid\partial f\mid^\frac{12}{5}+\frac{1}{r^6}\int_\Omega\eta^{s-1}e^{2nf}.\end{aligned}\end{equation}
Inserting these into (\ref{3-17}) and taking $\epsilon$ small we get
\begin{equation}\label{3-21}\begin{aligned}&\int_\Omega\eta^s\mathbf{Re}\left[ e^{2(n-1)f}\left(\mid\partial f\mid^2+e^{\frac{n}{n+m-1}(2+p)f}-2\frac{U_{\gamma\gamma}^{(k)}f_{k0}}{n}\right)\left(2U_{\beta\alpha}^{(k)}f_{\overline{\beta}}\right)\right]_\alpha f_{k0}\\&\lesssim\frac{1}{r^{4}}{\int}_{\Omega}\eta^{s-1}e^{2nf}e^{\frac{n}{n+m-1}(2+p)f}+\frac{1}{r^{6}}{\int}_{\Omega}\eta^{s-1}e^{2nf}.\end{aligned}\end{equation}
Similarly, other items of (\ref{3-10}) can be handled.Combining these with (\ref{3-9}) we can get
\begin{equation}\label{3-22}\begin{aligned}\int_\Omega\eta^s\mathcal{M}&\lesssim\frac{1}{r^2}\int_\Omega\eta^{s-1}e^{2nf}e^{\frac{2n}{n+m-1}(2+p)f}\\&\quad+\frac{1}{r^{4}}{\int}_{\Omega}\eta^{s-1}e^{2nf}e^{\frac{n}{n+m-1}(2+p)f}+\frac{1}{r^6}\int_\Omega\eta^{s-1}e^{2nf}.\end{aligned}\end{equation}

 To proceed, we need the following lemma.
\begin{lemma}\label{lem-1}
\begin{equation}\label{3-23}
\begin{aligned}&\int_{\Omega}\eta^{s}e^{\frac{3n}{n+m-1}(2+p)f+(2n-2)f}\\&\lesssim\int_{\Omega}\eta^{s}\mid\partial f\mid^{2}e^{\frac{2n}{n+m-1}(2+p)f+(2n-2)f}+\frac{1}{r}\int_{\Omega}\eta^{s-1}e^{\frac{2n}{n+m-1}(2+p)f+2nf}.\end{aligned}
\end{equation}
\end{lemma}
\vspace{10pt}
 \begin{proof}
 Multiply both sides of the equation (\ref{2.3}) by $-\eta^se^{\frac{2n}{n+m-1}(2+p)f+(2n-2)f}$
 and integrate over  $\Omega$ we have
\begin{equation}\label{3-24}\begin{aligned}&n\int_\Omega\eta^s\left(\mid\partial f\mid^2+e^{\frac{n}{n+m-1}(2+p)f}-\frac{U_{\alpha\alpha}^{(k)}f_{k0}}{n}\right)e^{\frac{2n}{n+m-1}(2+p)f+(2n-2)f}\\&=n\int_\Omega\eta^sge^{\frac{2n}{n+m-1}(2+p)f+(2n-2)f}\\&=-\int_{\Omega}\eta^sf_{\alpha\overline{\alpha}}e^{\frac{2n}{n+m-1}(2+p)f+(2n-2)f}\\&=\left[\frac{2n}{n+m-1}(2+p)+2n-2\right]\int_\Omega\eta^s\mid\partial f\mid^2e^{\frac{2n}{n+m-1}(2+p)f+(2n-2)f}\\&+s\int_{\Omega}\eta^{s-1}f_{\alpha}\eta_{\overline{\alpha}}e^{\frac{2n}{n+m-1}(2+p)f+\left(2n-2\right)f}.\end{aligned}\end{equation}
Since $U^{(k)}$ is a family of skew-Hermitian matrix, $U_{\alpha\alpha}^{(k)}$ is a pure imaginary number, then
\begin{equation}\label{3-25}\begin{aligned}&\int_\Omega\eta^se^{\frac{3n}{n+m-1}(2+p)f+\left(2n-2\right)f}\\&\lesssim\int_{\Omega}\eta^{s}\mid\partial f\mid^{2}e^{\frac{2n}{n+m-1}(2+p)f+\left(2n-2\right)f}\\&+\frac{1}{r}\int_{\Omega}\eta^{s-1}|\partial f|e^{\frac{2n}{n+m-1}(2+p)f+\left(2n-2\right)f}\\&\lesssim\int_{\Omega}\eta^{s}e^{\frac{2n}{n+m-1}(2+p)f+\left(2n-2\right)f}\\&+\frac{1}{r}\int_{\Omega}\eta^{s-1}e^{\frac{2n}{n+m-1}(2+p)f+\left(2n-2\right)f}.\end{aligned}\end{equation}
This is (\ref{3-23}). This ends the proof of Lemma \ref{lem-1}
\end{proof}

Combining this with (\ref{3-22}) and using Young's inequality we have
\begin{equation}\label{3-26}\begin{aligned}&\int_\Omega\eta^se^{\frac{3n}{n+m-1}(2+p)f+(2n-2)f}\\&\lesssim\frac{1}{r^{2}}{\int}_{\Omega}\eta^{s-1}e^{2nf}e^{\frac{2n}{n+m-1}(2+p)f}\\&\quad+\frac{1}{r^{4}}{\int}_{\Omega}\eta^{s-1}e^{2nf}e^{\frac{n}{n+m-1}(2+p)f}+\frac{1}{r^{6}}{\int}_{\Omega}\eta^{s-1}e^{2nf}\\&\lesssim\epsilon\int_{\Omega}\eta^{s-1}e^{\frac{3n}{n+m-1}(2+p)f+2nf}+r^{-2\times\frac{\frac{3n}{n+m-1}(2+p)+2n}{\frac{n}{n+m-1}(2+p)}}\int_{\Omega}2\eta^{s-1},\end{aligned}\end{equation}
by choosing  $s>0$ big enough and $\epsilon$ small, we finally obtain
\begin{equation}\label{3-27}\begin{aligned}&\int_{\mathrm{B}_{r}(\xi_{0})}\eta^{s}e^{\frac{3n}{n+m-1}(2+p)f+(2n-2)f}\\&\lesssim_{r}^{\bigg(2n+2m-2\times\frac{\frac{3n}{n+m-1}(2+p)+2n}{\frac{n}{n+m-1}(2+p)}\bigg)}\\&=r^{2n+2m-2\times\left(3+\frac{2\left(n+m-1\right)}{2+p}\right)}\\&=r^{Q-2\times[\frac{3q-q^*}{q-1}+\frac{2}{(q-q^*)(n+m-1)+2}]}.\end{aligned}\end{equation}

The theorem \ref{Thm2} is proved.
\qed
\\

{\bf Acknowledgments.} The authors would like to express their heartfelt thanks to Professors X Zhang and X Yang for their helpful comments and suggestions.

{\bf Conflicts of Interest:} The authors declare no conflict of interest.

{\bf Data Availability Statements:} No Applicable.

\end{document}